\documentclass[12pt]{amsart}
\usepackage{latexsym}

\newcommand{\haf}{{\frac{1}{2}}}
\newcommand{\holo}{\mathcal{O}}
\newcommand{\R}{{\mathbb R}}
\newcommand{\C}{{\mathbb C}}

\newcommand{\N}{{\mathbb N}}
\newcommand{\newsection}[1]
{\subsection{#1}\setcounter{theorem}{0} \setcounter{equation}{0}
\par\noindent}

\newtheorem{theorem}{Theorem}

\newtheorem{lemma}[theorem]{Lemma}
\newtheorem{corr}[theorem]{Corollary}

\newtheorem{proposition}[theorem]{Proposition}
\newtheorem{deff}[theorem]{Definition}
\newtheorem{example}[theorem]{Example}
\newtheorem{remark}[theorem]{Remark}

\newcommand{\eprop}{\end{proposition}}

\newcommand{\Sum}{\displaystyle\sum}
\newcommand{\Prod}{\displaystyle\prod}
\newcommand{\fract}{\frac}

\usepackage{setspace}

\begin{document}
\thispagestyle{empty}

\noindent {\large {\bf The zeros of Gaussian random holomorphic
functions on $\C^n$, and hole probability.\hfill\\
by Scott Zrebiec }}\medskip

\bigskip\bigskip
{\bf Abstract} \par We consider a class of Gaussian random
holomorphic functions, whose expected zero set is uniformly
distributed over $\C^n $. This class is unique (up to
multiplication by a non zero holomorphic function), and is closely
related to a Gaussian field over a Hilbert space of holomorphic
functions on the reduced Heisenberg group. For a fixed random
function of this class, we show that the probability that there
are no zeros in a ball of large radius, is less than $e^{-c_1
r^{2n+2}}$, and is also greater than $e^{-c_2 r^{2n+2}}$. Enroute
to this result we also compute probability estimates for the event
that a random function's unintegrated counting function deviates
significantly from its mean.
\newsection{Introduction}\par
Random polynomials and random holomorphic functions are studied as
a way to gain insight into difficult problems such as string
theory and analytic number theory. A particularly interesting case
 of random holomorphic functions is when the random functions can be defined so that they
 are invariant with respect to the natural isometries of the space in
question. The class of functions that we will study are the unique
Gaussian random holomorphic functions, up to multiplication by a
nonzero holomorphic function, whose expected zero set is uniformly
distributed on $\C^n$. 
For this class of random holomorphic functions we will determine
the expected value of the unintegrated counting function for a
ball of large radius and the chance that there are no zeros
present. This pathological event is what is called the hole
probability of a random function. In doing this we generalize a
result of Sodin and Tsirelison, to n dimensions, in order to give
the first nontrivial example, where the hole probability is
computed in more than 1 complex variable.
\par

The topic of random holomorphic functions is an old one which has
many results from the first half of the twentieth century, and is
recently experiencing a second renaissance. In particular Kac
determined a formula for the expected distribution of zeros of
real polynomials in a certain case, \cite{Kac43}. This work was
generalized throughout the years, and a terse geometric proof, and
some consequences are presented by Edelman and Kostlan,
\cite{EdelmanKostlan}. An excellent reference for other results
regarding the general properties of random functions is Kahane's
text, \cite{Kahane85}. One series of papers, by Offord, is
particularly relevant to questions involving the hole probability
of random holomorphic functions and the distribution of values of
random holomorphic functions, \cite{Offord67}, \cite{Offord72},
although neither is specifically used in this paper. Recently,
there has been a flurry of interest in the zero sets of random
polynomials and holomorphic functions which are much more natural
objects than they may initially appear. For example Bleher,
Shiffman and Zelditch show that for any positive line bundle over
a compact complex manifold, the random holomorphic sections to
$L^N$ (defined intrinsically) have universal high N correlation
functions, \cite{BleherShiffmanZelditchUniv}.
\par
In addition to a plethora of results describing the typical
behavior, there have also been several results in 1 (real or
complex) dimension for Gaussian random holomorphic functions where
the hole probability has been determined. For a specific class of
real Gaussian polynomials of even degree 2n, Dembo, Poonen, Shao
and Zeitouni have shown that for the event where there are no real
zeros, $E_n$, the $\lim_{n\rightarrow \infty}\frac{Prob
(E_n)}{\log (n)} n^{-b}= -b, \ b\in [0.4, 2]$,
\cite{DemboPoonenShaoZeitouni02}. Hole probability for the complex
zeros of a Gaussian random holomorphic function is a quite
different problem. Let $Hole_r=\{f$, in a class of holomorphic
functions, such that $\forall z\in B(0,r), \ f(z)\neq 0\}$. For
the complex zeros in one complex dimension, there is a general
upper bound for the hole probability: $Prob (Hole_r)\leq e^{-c \mu
(B(0,r))}, \ \mu(z)= E[Z_{\psi_\omega}]$ as in theorem
\ref{Onemoment}, \cite{Sodin00}. In one case this estimate was
shown by Peres and Virag to be sharp: $Prob(Hole_r) = e^{-\frac{
\mu(B(0,r))}{24}+ o(\mu(B(0,r)))}$, \cite{PeresVirag04}. These
last two results on hole probability might suggest that when the
random holomorphic functions are invariant with respect to the
local isometries, thus ensuring that $E[Z_\omega]$ is uniformly
distributed on the manifold, the rate of decay of the hole
probability would be the same as that which would be arrived at if
the zeros where distributed according to a Poisson process.
However, as the zeros repel in 1 dimension, \cite{Hannay96}, one
might expect there to be a quicker decay for hole probability of a
random holomorphic function. This is the case for random
holomorphic functions whose expected zero set is uniformly
distributed on $\C^1$, \cite{SodinTsirelison03} :
 $$Prob(Hole_r)\leq e^{-c_1 r^4}= e^{-c
\mu(B(0,r))^2}, \ and \ Prob(Hole_r)\geq e^{-c_2 r^4}= e^{-c
\mu(B(0,r))^2}$$ \par

The random holomorphic functions that we will study, can be
written as $$\psi_\omega(z_1, \ z_2, \ldots, \ z_n)= \Sum_{j\in
\N^n}\omega_j \frac{z_1^{j_1} \cdot z_2^{j_2} \cdot \ldots \cdot
z_n^{j_n}}{\sqrt{j_1 \cdot j_2 \cdot \ldots \cdot j_n}}=
\Sum_{j\in \N^n} \omega_j \frac{z^j}{\sqrt{j!}}
$$ where $\omega_j$ are independent identically distributed standard complex gaussian random variables,
and a.s. are holomorphic on $\C^n$.
The second form is just the standard multi-index notation, and
will frequently be used from here on out. Random holomorphic
functions of this form are a natural link between Hilbert spaces
of holomorphic functions on the reduced Heisenberg group and a
similar Gaussian Hilbert Space. Further, these random functions
will be the unique class (up to multiplication by a nonzero entire
function) whose expected distribution of the zero set is:
$$E[Z_\omega]= \frac{i}{2\pi}(dz_1\wedge d\overline{z}_1 +
dz_2\wedge d\overline{z}_2+\ldots +dz_n\wedge d\overline{z}_n) .
$$\par

The two main results of this paper are: \\

\begin{theorem}\label{Main}
If
$$\psi_\omega (z_1, z_2 \ldots, z_n)= \Sum_j \omega_j
\frac{z_1^{j_1} z_2^{j_2} \ldots z_n^{j_n}}{\sqrt{j_1! \cdot
j_n!}},$$ where $\omega_j$ are independent identically distributed
complex Gaussian random variables,
\newline then for all $\delta
> 0,$ there exists $c_{3,\delta}>0 \ and \ R_{n,\delta}$ such that for all $r> R_{n,\delta}$
$$Prob\left(\left\{\omega: \left|n_{\psi_\omega}(r) -\frac{1}{2}r^2
\right| \geq \delta r^2 \right\} \right) \leq e^{-c_{3,\delta}
r^{2n+2}}$$ where $n_{\psi_\omega}(r)$ is the unintegrated
counting function for $\psi_\omega$. \end{theorem}

\begin{theorem}\label{Hole probability}\par  If
$$Hole_r=\{ \omega:\forall z \in B(0,r), \
\psi_\omega (z)\neq 0 \},$$ then there exists $R_n, \ c_1, \ and \
c_2
>0$ such that for all  $r>R_n$
$$ e^{-c_2 r^{2n+2}} \leq Prob (Hole_r)\leq e^{-c_1 r^{2n+2}}$$
\end{theorem}



The proof of Theorems \ref{Main} and \ref{Hole probability}, will
use techniques from probability theory, several complex variables
and an invariance rule for Gaussian random holomorphic functions
which is derived from isometries of the reduced Heisenberg group.
These results, using the mainly same techniques, were already
proven in the case where n=1 by Sodin and Tsirelison,
\cite{SodinTsirelison03}.\par {\bf Acknowledgement:} I would like
to thank Bernie Shiffman for many useful discussions.


\newsection{The link between random holomorphic functions, Gaussian Hilbert spaces and the reduced Heisenberg group}\label{section:Background}
\par

To develop the notion of a random holomorphic function on $\C^n$
we will need a way to place a probability measure on a space of
holomorphic functions on $\C^n$. The definition we will use is
that a random holomorphic function is a representative of a
Gaussian field between two Hilbert spaces on the reduced
Heisenberg group. Through this definition we will prove the
crucial Lemma \ref{Invariance} which gives a nice law to determine
how random holomorphic functions behave under translation.
Additionally, this definition is equivalent to defining a random
holomorphic function as $\psi_\omega(z)= \Sum_{j\in \N^n} \omega_j
\left(\frac{z_1^{j_1} z_2^{j_2} \ldots z_n^{j_n}}{\sqrt{j_1! j_2!
\ldots j_n! }}\right)$, where $\omega_j$ are independent
identically distributed standard complex Gaussian random
variables.


We will start with the concept of a Gaussian Hilbert Space, as
presented by Janson, \cite{Janson97},
\begin{deff}
 A Gaussian Linear space, $G$, is a
linear space of random variables, defined on a probability space
$(\Omega, d \nu )$, such that each variable in the space is
Gaussian random variable.
\end{deff}\par

\begin{deff}
A Gaussian Hilbert Space, G is a Gaussian linear space that is
complete with respect to the $L^2(\Omega, d\nu)$ norm.
\end{deff}

We will shortly apply these definitions to a Hilbert space of
CR-holomorphic functions on the reduced Heisenberg group. The
Heisenberg Group, as a manifold is nothing other than
$\C^{n}\times \R$, and the reduced Heisenberg group is the circle
bundle:
 $X=H^n_{red}= \left\{(z, \alpha), z\in \C^n, \alpha \in \C, \ |\alpha|= e^{\frac{-|z|^2}{2}}
\right\}$.
\par Consider holomorphic functions of the Heisenberg group, which
are linear with respect to the $n+1^{st}$ variable. The
restriction of these functions define functions on X. For
functions on X there is the following inner product:
\par\begin{tabular}{rl} $\displaystyle{(F,G)= \int_X F \overline{G}}$ & $\displaystyle{= \frac{1}{\pi^n} \int_X f(z) \overline{g(z)} |\alpha|^{2} d\theta (\alpha) dm(z) }$
\\ & $\displaystyle{= \frac{1}{\pi^n} \int_{\C^n} f(z) \overline{g(z)} e^{-|z|^2} dm(z) }$
\end{tabular}\\ Here $dm$ is Lebesque measure. With
respect to this inner product, \par $H_X=\{ F \in
\mathcal{O}(\C^{n+1}), \ F(z,\alpha) = f(z) \alpha, \ f \in
\holo(\C^n) \}$ is a Hilbert Space, and $H_X \cong H^2(\C^n, \
e^{-|z|^2} dm) $, as Hilbert Spaces.

\begin{proposition} For $H_{X}, \ \left\{ \frac{z_1^{j_1}\cdot
z_2^{j_2}\cdot \ldots z_n^{j_n} }{\sqrt{j_1! \cdot \ldots \cdot
j_n! }} \alpha \right\}_{j\in \N^n}= \{ \psi_j(z) \alpha\}_{\N^n}
$ is an orthonormal basis.
\end{proposition}

The proof of this proposition is a straight forward
computation.\par

The isometries of the reduced Heisenberg group will play a crucial
role in my computation of the hole probability. These isometries
are of the form:
$$ \tau_{(\nu, \alpha)}: H^m_{red} \rightarrow H^m_{red} $$
$$ \tau_{(\nu, \alpha)}: (\zeta, \beta)  \rightarrow  (\zeta + \nu, e^{-\zeta \overline{\nu} } \alpha \beta) $$

\par

The inner product on $X$ is invariant with respect to the
Heisenberg group law:\\ \begin{tabular}{rl} $\tau^*\langle
F,G\rangle $&$= \frac{1}{\pi^n} \int \int f(\zeta+ \nu) \overline
{g(\zeta+ \nu)} |\beta \alpha e^{-\zeta \overline \nu}|^{2}
d\theta \ dm(\zeta) $\\&$= \frac{1}{\pi^n} \int f(\zeta+ \nu)
\overline{g(\zeta+ \nu)} e^{-|(\zeta + \nu)|^2} dm(\zeta + \nu)
$\\& $= \frac{1}{\pi^n} \int f(\zeta)
\overline{g(\zeta)} e^{-|\zeta|^2 } dm(\zeta) $\\
\end{tabular}\\

As such for $\alpha= e^{-\frac{|z|^2}{2}}$,
$$\tau^* (\alpha \psi_j(z))= (\alpha e^{-\haf|\zeta|^2- z \overline{\zeta}}
\psi_j(z+\zeta))= e^{-\haf |z+\zeta|^2-i\cdot Im (z
\overline{\zeta})} \psi_j(z+\zeta)$$ and the collection of these,
$\{e^{-\haf |z+\zeta|^2- i\cdot Im (z \overline{\zeta})}
\psi_j(z+\zeta)\} $, is another orthonormal basis for $H_X$, as
the inner product is invariant with respect to the group law.
\par

\begin{example} (Gaussian Hilbert Spaces)\par Let $G_{H_{X}}'= Closure(Span\left( \{\omega_j
\psi_j(z) \alpha\}_{\N^n}\right))$, where the closure is taken
with respect to the norm $ E[ (\| \cdot \|_{H_X})^2]^\haf $ and
where $ \omega_j$ are independent identically distributed standard
complex Gaussian random variables. $G_{H_{X}}'$ is not a Gaussian
Hilbert space but is isometric to $G_{H_{X}}= Closure(Span\left(
\{\omega_j \}_{\N^n}\right))$, which is.\par
\begin{tabular}{rl}Of course,&$\ \ \ \ \ \ \ \ H_X \ \ \ \ \ \rightarrow \
\ \ \ \ G_{H_X}'\ \ \ \ \ \ \ \ \ \rightarrow \ \ \ G_{H_X} $\\&
$\ \Sum a_j \psi_j(z)\alpha \longmapsto \Sum a_j \omega_j
\psi_j(z)\alpha
 \longmapsto \Sum a_j \omega_j$
\end{tabular}\\ are isometries.\end{example}
\par $G_{H_X}'$ is in many ways more natural then $G_{H_X}$,
and is closely related to random holomorphic functions.

\begin{deff}
A Gaussian field is a linear isometry $L: H_X \rightarrow G_{H_X}
$. As such, for all $\displaystyle{f\in H_X, L[f]= X_f, a \
standard \ complex \ Gaussian \ random \ variable}$\\
$\displaystyle{ \ with \ Var = \| f \|_{H_X}^2}$.\end{deff}

\begin{deff}
 A Gaussian random function, is a representative for a
Gaussian field $L$. In other words,
$$if \  f\in H_X, L[f]=\langle \varphi_\omega , f \rangle_{H_X}$$
\end{deff}


\begin{remark} (Random holomorphic functions on $\C^n$ or (equivalently) Gaussian fields and functions between $H_{X}$ and $G_{H_{X}}$)
\par Let
$\psi_{\omega}(z) = \Sum_{j\in \N} \omega_j \psi_j(z)= \Sum_{j\in
\N} \omega_j \frac{z^j}{\sqrt{j!}}$ , $\omega_j$ independent
identically distributed standard complex Gaussian random
variables. This will shortly be shown to a.s. be a holomorphic
function on $\C^n$ by Theorem \ref{Convergence}.\par Note that for
$\alpha f(z)\in H_{X}, \ f(z)= \Sum a_j \psi_j(z), \ \{a_j\}\in
\ell^2$,\newline $\langle \alpha \psi_\omega(z), \alpha f(z)
\rangle= \Sum \overline{a}_j \omega_j $, which is a complex
Gaussian random variable with variance $\Sum |a_j|^2=
\|f\|_{H_{\C^n}} $, hence $\alpha \psi_\omega(z)$ is a random
CR-holomorphic function on $X$.\par The variable $\alpha$ will be
useful when we change bases. This occurs when we look at how
random functions behave with respect to translation (Lemma
\ref{Invariance}). Abusing notation we will frequently drop it and
we will call $\psi_\omega(z)$ a random holomorphic function on
$\C^n$.
\par
\end{remark}

 There is a simple condition
for when a function of the form $\Sum \omega_j \psi_j(z)$ is a
holomorphic condition, where $ \omega_j$ are independent
identically distributed complex Gaussian Random variables.
\par

\begin{theorem}\label{Convergence}
Let $\{\omega_j\}_{j\in \N}$ be a sequence of independent
identically distributed, standard complex Gaussian random
variables. If for $j\in \N, \newline \psi_j(z)\in \mathcal{O} (
\Omega)$, and for all $\ compact \ K\subset \Omega,$
$\displaystyle{\Sum_{j\in \N} \max_{z \in K} |\psi_j(z)|^2<
\infty}$\par then for a.a.-$\omega, \ \Sum_{j\in \N}\omega_j
\psi_j(z)$
 defines a holomorphic
function on $\Omega$.
\end{theorem}

This theorem can be proved easily by adapting a similar proof of
convergence of "random sums" from \cite{Kahane85}.

\begin{theorem}\label{GRHF defnthm}
If $L$ is a Gaussian field, $L: H_X \rightarrow G_{H_X}$, and $\{
\phi_j\}$ is an orthonormal basis for $H_X$
\par then $L$ can be written as: $L[\cdot]= \langle \phi_X,
\cdot \rangle $, where $\phi_X(z)= \sum X_i \phi_i$, and $\{X_i\}
$ is a set of independent identically distributed standard
Gaussian random variables.
\end{theorem}

\begin{proof} Let $L$ be a Gaussian random functional.\\
Let $X_1 = L[\phi_1], \ X_2 = L[\phi_2], \ \ldots , \ X_j =
L[\phi_j], \ldots $ \par We must only show that $X_j$ are
independent identically distributed Gaussian Random variables,
hence it suffices to prove independence as by the definition of
Gaussian random field, $X_i, \ X_j$ are jointly normal, as
$L\left[\Sum a_j \psi_j\right]= \Sum a_j X_j $ is normal.\\  For
$i\neq j$:\par \begin{tabular}{rl} $2$ &$=E[|X_i +X_j|^2]=
E[|X_i|^2]+ E[|X_j|^2] + E[X_i \overline{X}_j]+ E[X_j
\overline{X}_i] $\\  &$= 2 + E[X_i \overline{X}_j]+ E[X_j
\overline{X}_i]$\end{tabular}
\\  Hence
$Re(E[X_i \overline{X}_j])=0=Im(E[X_i \overline{X}_j])$,\par The
result then follows.\\
\end{proof}

\newsection{Common Results}\par

Let us briefly review properties of the zeros of random
holomorphic functions. An elementary way to view the zeros of a
holomorphic function is as a set: $Z_f= f^{-1}(\{0\})$, but this
will be insufficient for my purpose, and we will instead view it
as a (1,1) current. For $M^n$ an n dimensional manifold, and $f\in
\holo(M), \ f: M^n \rightarrow \C, \ f^{-1}(\{0\})$ is a divisor.
Hence the regular points of $Z_f$ are a manifold, and by taking
restriction we identify forms in $D^{(n-1,n-1)}_M$ with ones in
$D^{(n-1,n-1)}_{Z_{f,reg}}$. As $Z_{f, reg}$ is an n-1 complex
manifold, $\displaystyle{\int_{Z_{f, reg}} }$ is a (1,1) current
on M, which we will denote $Z_f$ (abusing notation). As the
singularities occur in real codimension 2. $Z_f= Z_{f,reg}$, and
in general: \label{Poincare-Lelong formula}if $f\in\holo(M^n)$, M
an n complex manifold, then $Z_f= \frac{i}{2 \pi}\partial
\overline{\partial} \log |f|^2$, as $(1,1)$ currents on M.\par

Before we classify the atypical hole probability, we shall first
describe the expected behavior. Many various forms of the
following theorem
 have been proven, \cite{EdelmanKostlan}, \cite{Kac43} and \cite{Sodin00}. For my purposes it is important that the proof
 is valid in n-dimensions, and for infinite sums. Many of the proofs resemble this one. After a
 conversation with Steve Zelditch, I was able to simplify a
 previously complicated argument into the current form. This simplification is already known to other
 researchers including Mikhail Sodin. \par For the following theorem let $\psi_j: \Omega
\rightarrow \C$, $j \in \Lambda, \ \Lambda = \{0, 1, 2, \ldots, n
\}$ or $\Lambda=\N$, be a sequence of holomorphic functions on a
domain of an n manifold to $\C$. \par

\begin{theorem}\label{Onemoment}
If $ E[|\psi_\omega|^2]=\Sum |\psi_j(z)|^2 \ converges \ locally \ uniformly \ in \  \Omega$\\
then $E[Z_\omega] = \fract{i}{2\pi}\partial\overline{\partial}
\log ||\psi(z)||_{\ell^2}^2$
\begin{proof}
\par
Let $\beta \ \in \ D^{n-1,n-1}(\Omega)$\par To simplify the
notation, let $\beta=\phi \ dz_2 \wedge d\overline{z}_2 \wedge
\ldots \wedge dz_n\wedge d\overline{z}_n$.\par

\begin{tabular}{ccl}$ \langle Z_{\psi_\omega},\beta  \rangle $ & $ = $& $  \langle \frac{i}{2\pi}
\partial\overline{\partial} \log (|\psi_\omega(z)|^2) ,  \beta
\rangle$\\
& $=$&$ \langle \frac{i}{2\pi} \log
\left(|\psi_\omega(z)|^2\right) ,
\partial
 \overline{\partial} \beta \rangle$\\

& $ =$ & $ \langle \frac{i}{2\pi} \left(\log
(||\psi(z)||_{\ell^2}^2) + \log
\left(\frac{|\psi_\omega(z)|^2}{||\psi(z)||_{\ell^2}}
\right)\right) ,
\partial \overline{\partial} \beta \rangle$\end{tabular} \par
Taking the expectation of both sides we compute:\par
\begin{tabular}{cl}
 $E[\langle Z_{\psi_\omega},\beta  \rangle] = \frac{1}{2\pi}$&$  \int_{\omega} \int_{z \in \Omega}
\log (||\psi(z)||_{\ell^2}^2) \frac{\partial^2 \phi} {\partial z_1
\partial \overline {z}_1} \ dm(z) \ d\nu(\omega)$\\ & $+ \frac{1}{2\pi}  \int_{\omega} \int_{z \in \Omega} \log
\left(\frac{|\psi_\omega(z)|^2}{||\psi(z)||_{\ell^2}^2}\right)
\frac{\partial^2 \phi} {\partial z_1 \partial \overline {z}_1} \
dm(z) \ d\nu(\omega)$\end{tabular}\par The first term is the
desired result (which by assumption is integrable and finite),
while the second term will turn out to be zero. We first must
establish that it is in fact integrable:
\par
\begin{tabular}{cl}
$ \int_{z \in \Omega} \int_{\omega} |\log
\left(\frac{|\psi_\omega(z)|^2}{||\psi(z)||_{\ell^2}^2} \right)$&$
\frac{\partial^2 \phi}{\partial z_1 \partial \overline {z}_1}
d\nu(\omega) | \ dm(z)$\\&$\leq c \int_{z \in K} \int_{\omega}
\left| \log
\left(\frac{|\psi_\omega(z)|^2}{||\psi(z)||_{\ell^2}^2}\right) \
d\nu(\omega)\right| \ dm(z)$\\& $= c \int_{z \in K} \int_{\omega}
\left| \log (|\omega'|^2) \right| \ d\nu(\omega') \
dm(z)$\end{tabular}\par where $\omega '$ is a standard centered
Gaussian ($\forall z$), thusly proving integrability as: $$
\int_{z \in \Omega} \int_{\omega} \left| \log
(\frac{|\psi_\omega(z)|^2}{||\psi(z)||_{\ell^2}^2}) \right| \leq C
\int_{|x|<1} | \log(x) | dm(x)+ c \int_{|x|>1} |xe^{-x^2}|
dm(x)\leq c$$
\\
Finally, $$\int_\Omega \beta \wedge E[Z_{\psi_\omega}]= \frac{i}{2
\pi} \int_\Omega
\partial \overline{\partial}\beta \log(\|\psi(z)\|_{\ell^2}^2) \ +
\ \int_\Omega C \frac {\partial^2 \phi}{\partial z_i \partial
\overline{z}_j} dm(z)$$
$$=\frac{i}{2 \pi} \int_{\Omega}\partial \overline{\partial} \beta \wedge \log(||\psi(z)||_{\ell^2}^2) $$
\end{proof}
\end{theorem}

\begin{corr}
For $\psi_\omega $ a random holomorphic function on $\C^n$,
$$E[Z_{\psi_\omega}] = \frac{i}{2\pi}(dz_1\wedge d\overline{z}_1 +
dz_2\wedge d\overline{z}_2+\ldots +dz_n\wedge d\overline{z}_n)$$
\end{corr}

In Theorem \ref{Onemoment}, we proved that the expected zero set
is determined by the variance of a of the random function when
evaluated at a point. More can be said:

\begin{theorem}\label{uniquenessOfIsoInvarGrhf}
For gaussian analytic functions the expected zero set determines
the process uniquely (up to multiplication by nonzero holomorphic
functions) on a simply connected domain.
\end{theorem}
This theorem is proven in one dimension by Sodin \cite{Sodin00},
and the same proof works in n-dimensions.
\newsection{Invariance of Gaussian random functions with respect to the isometries of the reduced Heisenberg
group}\label{InvarianceSection}\par The invariance property of the
random function in question with respect to the reduced Heisenberg
group's isometries plays a central role in proving that:
$$\displaystyle{{\max_{z\in
\partial B(0,r)}\left(\log(|\psi_\omega(z)| ) - \haf |z|^2\right)} =
{\max_{z\in \partial B( \zeta ,r)} \left(\log(|\psi_\omega(z) |) -
\haf |z|^2\right)} }$$ This invariance property, which was known
and used in the 1 dimensional case, and makes sense (from the view
that $\forall (z,\alpha)\in X, \ \alpha \psi_\omega(z)$) defines a
standard complex Gaussian random variable for any fixed $z$.
Apparently, however, there was no proof in the literature until
recently, \cite{Sodin05}. I also independently came up with this
same result by using the properties of the Heisenberg group, and
this is presented here.\par

This next result is that a random holomorphic function is well
defined independent of basis chosen, and also will be shortly
restated in order to give an important translation law for random
holomorphic functions of $\C^n$.
\begin{lemma}\label{IndepOfBasis}
If $\{\phi_j \}_{j \in
\Lambda} $ is an orthonormal basis for $H_X$ 
\par then there exists $\{\omega_j'\}_{j\in \Lambda}$ independent identically distributed standard complex
Gaussian random variables such that $for \ all \ (z,\alpha) \in H^n_{red},$  $$\alpha \psi_\omega(z) =
\phi_{\omega'}(z,\alpha), \  a.s.$$ where
$\phi_{\omega'}(z,\alpha)= \sum \omega_j' \phi_j(z,\alpha)$.
\end{lemma}

\begin{proof}
For all $j\in \N$, let $\omega_j' = \langle \alpha \psi_\omega(z),
\phi_j (z,\alpha) \rangle_{H_X}$, which is a standard complex
Gaussian random variable by Theorem \ref{GRHF defnthm}. Further,
for $j\neq k, \ \omega'_j$ and $\omega'_k$ are independent.\par
Let $f\in H_X \Rightarrow f = \sum a_j \phi_j, \ \{a_j\}_{j\in \N}
\in \ell^2$
\\ We now demonstrate that $\phi_{\omega'} = \psi_\omega$ as a
Gaussian field:\par
\begin{tabular}{rl} $\langle \phi_{\omega'}, f\rangle $&$=\sum \overline{a}_j
\omega_j'$\\
$\langle \psi_{\omega}, f\rangle $&$= \langle \psi_{\omega}, \sum a_j \phi_j \rangle$\\
&$= \sum \overline{a}_j \langle \psi_{\omega}, \phi_j \rangle $\\
&$= \sum \overline{a}_j \omega_j' $\\
\end{tabular}\par As $H_X= H_X^*$, for evaluation maps $eva_{(z_0,\alpha_0)}= \langle \sum b_n \phi_j, \cdot \rangle,$\newline $ \sum b_n \phi_j\in H_X$ and therefore by the above work:\par
$\langle \phi_{\omega'}, \sum \overline{b}_j \phi_j \rangle = \sum
\omega_j' b_j = \sum \omega'_j \phi_j(z_0,\alpha_0)=  \sum
\omega_j \alpha_0 \psi_j(z_0)$
\end{proof}


\begin{deff}
A Gaussian random function is invariant with respect to $\tau$ if
both $\langle\cdot, \psi_\omega \rangle$ and $\langle\cdot, \tau^*
\psi_\omega \rangle$ induce the same Gaussian
field.\\
\end{deff}

A random CR-holomorphic function on X, will be invariant with
respect to isometries of X. These will in turn be important for
random holomorphic functions on $ \C^n$, which is illustrated in
the following simple but important lemma.\par

Let $\tau(z, \alpha)= (z+ \zeta, e^{-z \overline{\zeta}} \alpha
\beta ), \ |\beta|=e^{-\frac{|\zeta|^2}{2}}$. Recall that $\tau$
is an isometry of $H_X$.
\par

\begin{lemma}\label{Invariance}
For all $z \in \C^n,$  there exists $\ \omega_j'$ independent
identically distributed standard complex Gaussian random
variables, such  that
 $$ \ e^{-\haf|z|^2} \psi_\omega(z) = e^{-\haf |z+ \zeta|^2 - i \cdot Im(z\overline{\zeta})}\psi_{\omega'}(z+\zeta)$$\\
\end{lemma}

\begin{proof} Here, $\zeta$ is any fixed complex number and we set $\beta=e^{-|\zeta|^2}$.
Both $\{\alpha \psi_j(z)\} $ and $\{\tau^*( \alpha \psi_j(z))\}$ are orthonormal bases of $H_X$,
and they therefore induce the same Gaussian random function, as
these are well defined independent of basis by Lemma
\ref{IndepOfBasis}. Hence\par
\begin{tabular}{rrl}
$\displaystyle{\alpha \psi_\omega(z)}$&$\displaystyle{=}$&$\alpha
\beta e^{-z\overline{\zeta}}\psi_{\omega'}(z+\zeta)$\\
$\displaystyle{e^{-|z|^2}
\psi_\omega(z)}$&$\displaystyle{=}$&$\displaystyle{e^{-\haf
(|z|^2+ |\zeta|^2 + 2 z \overline{\zeta})}
 \psi_{\omega'}(z+\zeta)}$\\
&$\displaystyle{=}$&$\displaystyle{e^{-\haf|z + \zeta|^2 - i \cdot
Im( z \overline{\zeta})}
 \psi_{\omega'}(z+\zeta)}$\\
\end{tabular}

\end{proof}

\begin{corr}
The random variable: $\displaystyle{{\max_{z\in
\partial B(0,r)} \left(\log(|\psi_\omega(z) |) - \haf |z|^2 \right) }  }$ is invariant with
respect to $\tau^*$. In other words
$$\displaystyle{{\max_{z\in
\partial B(0,r)}\left(\log( |\psi_\omega(z)| ) - \haf |z|^2\right)}  =
{\max_{z\in \partial B( \zeta ,r)} \left(\log( |\psi_\omega(z) |)
- \haf |z|^2\right)} }$$
\end{corr}

\begin{proof} This corollary just specializes the previous lemma
as,\\
\begin{tabular}{ll}
$\displaystyle{\max_{z\in \partial B(0,r)}} \left(\log |
\psi_\omega(z) | - \haf |z|^2\right)$ & $=
\displaystyle{\max_{z\in
\partial B(0,r)}} \left(\log| \psi_\omega(z) | - \log(e^{\haf |z|^2})\right)$\\
& $=\displaystyle{\max_{z\in
\partial B(0,r)}} \log(|\alpha
\psi_\omega(z) |), \ |\alpha|= e^{-\haf|z|^2}$ \\ &
$=\displaystyle{\max_{z\in
\partial B(0,r)}} \log(
\tau^*(|\alpha \psi_\omega(z) | ))$\\& $= \displaystyle{\max_{z\in
\partial B(0,r)}} \log( |\beta \psi_\omega(z+\zeta) | ), \ |\beta|= e^{-\haf |z+\zeta|^2}$ \\
& $=\displaystyle{\max_{z\in
\partial B(\zeta,r)}} \log( (|\beta' \psi_\omega(z) | )), \ |\beta'|= e^{-\haf |z|^2}$ \\
& $=\displaystyle{\max_{z\in
\partial B(\zeta,r)}} \left(\log |\psi_\omega(z)|- \haf |z|^2 | \right)$
\end{tabular}\\

\end{proof}


\newsection{An estimate for the growth rate of random holomorphic functions on the reduced Heisenberg group}\par
In this section we begin working towards my main results. Lemma
\ref{Anchises} is interesting in and of itself as it proves that
random functions for $H_X$ are of finite order 2, a.s. From hence
forth we will work with $\C^n$, for any one fixed n.
\begin{deff}
A family of events $\{E_r\}_{r\in \R^+}$, dependent on r, will be
called a small family of events if exists $R, \ and \ c>0, \ such
\ that \ $ for all $r>R, \ Prob(E_r)\leq e^{-cr^{2n+2}}$.

\end{deff}

We will be using properties of Gaussian random holomorphic
functions to deduce typical properties of functions, and the size
of the family of events where these typical properties will not
work will always be small.\par Let $\displaystyle{M_{r,
\omega}=\max_{\partial B(0,r)} \log |\psi_\omega(z)|}$
\par We will be able to compute this, adapting a strategy that
Sodin and Tsirelison, \cite{SodinTsirelison03}, used to solve the
analogous 1 dimension problem, by using the Cauchy Integral
Formula in conjunction with some elementary probability theory and
computations:

\begin{lemma}\label{Gauss}
Let $\omega$ be a standard complex Gaussian RV (mean 0,
Variance 1), with a probability distribution function $d\nu(\omega)$ \\
\begin{tabular}{cl}
then: & a-i) $\nu(\{ \omega : | \omega| \geq \lambda \}) =
e^{-\lambda^2}$\\
& a-ii) $\nu(\{ \omega : | \omega| \leq \lambda \}) =1-
e^{-\lambda^2}\in [\frac{\lambda^2}{2}, \lambda^2], if \lambda
\leq 1$
\end{tabular}\par
b) If $\{\omega_j \}_{j\in \N^n}$ is a set of of independent
identically distributed standard Gaussian random variables, then
$\nu(\{ \omega: |\omega_j| < (1+ \varepsilon)^{|j|}\})=c>0$.
\end{lemma}\par Here, and throughout this paper for $j\in \N^n, \ |j|:=\sum j_i $\par
Lemma \ref{Gauss}-a) is a straight forward computation using that
the probability distribution for a standard complex Gaussian.
\begin{proof} of b)
$\nu (\{\omega_j: |\omega_j| < (1+\varepsilon)^{|j|} \})=
1-e^{-(1+\varepsilon)^{2|j|}}$\\\begin{tabular}{rl} $\left(\nu(\{
\omega: |\omega_j| < (1+ \varepsilon)^{|j|}\})>0\right)$ & $\Leftrightarrow \Prod_{j\in \N^n, |j|=0, }^{|j|=\infty} 1-e^{-(1+\varepsilon)^{2|j|}} = c $ \\
&$\Leftrightarrow c \Sum |j|^n
|\log\left(1-e^{-(1+\varepsilon)^{2|j|}} \right)|<\infty
 $\end{tabular}\\ as there are about $c |j|^n, \ j\in \N^n$ with a fixed value of $|j|$.\par
$\forall |x|<1, \ \log(1-x)= - \int  \Sum (x)^m =\Sum_{m\geq 0}
\frac{x^{m+1}}{m+1}$\par Therefore, $\Sum |j|
\left|\log\left(1-e^{-(1+\varepsilon)^{2|j|}} \right)\right| \leq
c \Sum_m m^n e^{-(1+\varepsilon)^{2m}}< \infty $
\end{proof}

The following lemma is needed twice in this paper, including in
the proof of Lemma \ref{Anchises}.

\begin{lemma}\label{CalcII}
If $j \in \N^{+,n}$ then $\frac{|j|^{|j|}}{j^{j}} \leq n^{ |j|}$
\end{lemma}

\begin{proof}
Let $u_k = \frac{j_k}{|j|}\geq 0$, hence $\sum_{k=1}^n u_k =1$.\\
As $\sum u_k \delta_{u_k}$ is a probability measure:
\par \begin{tabular}{rl} $\displaystyle{ \sum_{j=1}^n u_k \log
\left(\frac{1}{u_k}\right)}$&$\displaystyle{\leq \log \sum
\frac{1}{u_k} u_k
}$, by Jensen's inequality.\\
&$\displaystyle{= \log (n)}$
\end{tabular}\\
\begin{tabular}{rl}Hence, $
n^{|j|}$&$\displaystyle{ \geq \prod_k ( u_k)^{-|j| u_k}}$ \\&
$=\frac{ |j|^{|j|}}{ j^j } $
\end{tabular}



\end{proof}

\begin{lemma}\label{Anchises} (Probabilistic Estimate on the Rate of growth of the maximum of a random function on $\C^n
$)\\
For all $\delta >0$,
$$E_{r,\delta}:= \left\{\omega : \left|\frac{log (M_{r,\omega})}{r^2}- \frac{1}{2} \right| \geq \delta \right\} \ is \ a \
small\ family \ of \ events$$

\end{lemma}

\begin{proof}
We will first prove that: $\nu(\{\omega : \frac{log
(M_{r,\omega})}{r^2} \geq \frac{1}{2}+ \delta \} \leq e ^ {-
c_{\delta, 1} r^{2n+2} }$ and we will prove this by specifying a
set of measure almost 1 where the max grows at the appropriate
rate.
\\
\begin{tabular}{ccl}
Let $\Omega_r$ be the event where: & $i) \ |\omega_j|\leq
e^{\frac{\delta r^2}{4}},$ & $ |j| \leq
2 e \cdot n \cdot r^2  $\\
& $ii) |\omega_{j}|\leq 2^{\frac{|j|}{2}}, $ & $ \ |j| > 2 e \cdot n \cdot r^2$\\
\end{tabular}
\\

\begin{tabular}{cl}$\nu(\Omega_r^c)$ & $\displaystyle{\leq \sum_{|j|\leq 2 e \cdot n r^2} \nu ( \{|\omega_j|> e^{\frac{\delta r^2}{4}} \})+ \sum_{|j| > 2 e \cdot n  r^2} \nu ( \{|\omega_j|> 2^{\frac{|j|}{2}} \}) } $ \\
& $\displaystyle{ \leq c_n r^{2n} e^{\left(-e^{\frac{\delta r^2}{2} }\right)}+\Sum_{|j|_> 2 e \cdot n  r^2} e^{-2^{|j|}}}$\\
& $\leq e^{-e^{cr^2}} + c e^{-2^{c r^2}}, \ \forall r>R_0  $ \\ & $ \leq e^{-e^{c r}}$\\
\end{tabular}
\par We now have that $\Omega_r^c$ is contained in a small family of
events (and in fact could make a stronger statement on the rate of
decay in terms of r). It now remains for me to show that $\forall
\omega\in \Omega_r, \ \frac{\log|M_{r,\omega}|}{r^2}\leq \haf +
\haf \delta.$
\par
$\forall \omega \in \Omega_r,$ we have that:\\
$M_{r,\omega} \leq \Sum_{|j| = 0}^{|j|\leq 4 e \cdot n (\frac{1}{2} r^2)} |\omega_j| \frac{ |z|^{j}}{\sqrt{j!}}+ \Sum_{|j| > 4 e \cdot n (\frac{1}{2}r^2)} |\omega_j| \frac{ |z|^{j} }{\sqrt{j!}} = \sum^1 + \sum^2$ \\
\vspace{.1in}\\Using the Cauchy-Schwartz inequality: \par
\begin{tabular}{cl}$\Sum^1 $
& $ \leq (e^{\frac{1}{4} \delta r^2} ) \sqrt{c (r^2)^n} \left(\Sum_j \frac{|z^{2 j}|}{j!}\right)^\frac{1}{2}$\\
& $\leq c_n e^{\frac{\delta r^2}{4}} r^n e^{\frac{1}{2}r^2} $ \\
&$ \leq
e^{(r^2)(\frac{1}{2}+\frac{1}{3}\delta)}, \ \forall r> R_{n,\delta}$.\\
\end{tabular}
\par
\begin{tabular}{cl}$ \Sum^2$
& $ \leq \Sum _{|j| > 4 e \cdot n r^2} (2)^{\frac{|j|}{2}} \frac{ |z^{j}| }{\sqrt{j!}}$ \\
& $ \leq \Sum _{|j| > 4 e \cdot n r^2} (2)^{\frac{|j|}{2}} \left(\frac{|j|}{4 e n}\right)^{\frac{|j|}{2}} \prod_k \left(\frac{e}{j_k}\right)^{\frac{j_k}{2}} $, by Sterling's Formula\\
& $\leq C $, by Lemma \ref{CalcII}.\end{tabular} \par Hence,
$\forall \omega \in \Omega_r, \ \log(M_{r,\omega})\leq
(\frac{1}{2}+\frac{1}{2}\delta)r^2 $\par

It now remains for me to show that:
$$\forall \delta< \Delta, \ \nu\left( \left\{\omega : \frac{log (M_{r,\omega})}{r^2} \leq \frac{1}{2}- \delta \right\}\right) \leq e ^ {- c_{\delta, 2} r^{2n+2} }$$
which we will do by using Cauchy's integral formula to transfer
information on $M_{r, \omega}$ to individual coefficients
$\omega_j$. It suffices to prove this result only for small
$\delta$ as $\delta< \delta' \Rightarrow E_{\delta', r}\subset
E_{\delta, r}$. The constant $\Delta$ can be explicitly
determined.\par It will be most convenient to prove this result
for the polydisk, where the Cauchy Integral Formula applies. The
notation for the polydisk is the standard one: $P(0,r):=\{z\in
\C^n: \ \forall i, \ |z_i|<r \} $\par Let
$\displaystyle{M_{r,\omega}'= \max_{z\in P(0,r)}
|\psi_\omega(z)}|$\par The corresponding claim for a poly disk is
that: $$\displaystyle{M_{r,\omega}'}\geq \frac{n}{2} r^2- \delta
r^2$$except for a small family of events.
\par We will now look at the probability of the event consisting of $\omega $
such that:
$$\log(M_{r,\omega}') \leq \left(\frac{n}{2}-\delta \right)r^2  $$
By Cauchy's Integral Formula: $\left|\frac{\partial^{ j } \psi_\omega}{\partial z^j}\right|(0) \leq j! M_{r,\omega}' r^{-|j|}$ \\
By direct computation using the definition of $ \psi_\omega(z) $
in terms of a power series:$$\left|\frac{\partial^{ j }
\psi_\omega}{\partial z^j}\right|(0)= |\omega_j| \sqrt{j!}
$$\par
Therefore: 
 $|\omega_j| \leq c M_{r,\omega}' \sqrt{j!} r^{-|j|},$
\\
and using Sterling's formula ($j! \approx \sqrt{2 \pi} \sqrt{j}
j^j e^{-j} $), we get that:\par $|\omega_j|\leq (2
\pi)^{\frac{n}{2}} (\prod_k j_k^{\frac{1}{4}}) e
^{{(\frac{n}{2}-\delta)r^2}+\sum \frac{{j_k}}{2} \log({j_k}) -(|j|) \log r- \frac{|j|}{2}} $, $\forall k , \ j_k\neq 0.$ \\
 \\The $(2\pi)^{\frac{n}{2}}j^{\frac{1}{4}}$ term will not matter in the end so we will focus instead on the exponent.
\\

\begin{tabular}{rl}$\displaystyle{A}$&$=\displaystyle{
\left(\frac{n}{2}-\delta\right)r^2-\frac{|j|}{2}+ \Sum_k
\left(\frac{j_k}{2}\log(j_k)\right) - (|j|) \log(r)}$
\\ &$\displaystyle{= \Sum_{k=1}^{k=n} \left(\frac{j_k}{2}\right) \left(\left(1-\frac{2\delta}{n}\right)\frac{r^2}{j_k}-1+
 \log(j_k)
- 2 \log(r)\right)}$\end{tabular}\\
Let $j_k = \gamma_k r^2$\par

\begin{tabular}{rl}$\displaystyle{ A}$ & $\displaystyle{ = \Sum_{k=1}^{k=n} \left(\frac{\gamma_k
r^2}{2}\right)
\left(\left(1-\frac{2\delta}{n}\right)\frac{1}{\gamma_k }-1+
 \log(\gamma_k)
\right)}$ \\ & $\displaystyle{= - \delta r^2+ n f(\gamma_k)
\frac{r^2}{2}}$, where $f(\gamma_k) = 1- \gamma_k + \gamma_k \log
(\gamma_k)$
\end{tabular}\\ \par $f(\gamma_k) = (1-\gamma_k)^2 - (1-\gamma_k)^3 +
o((1-\gamma)^4)$ near 1.\par

Hence $\exists \Delta$ such that $\forall \delta\leq \Delta$ if
$\gamma_k \in \left[1-\sqrt \frac{\delta}{n}, 1+\sqrt
\frac{\delta}{n} \right]$ then $A\leq \frac{-\delta r^2}{2} $ \par

Therefore for $j$ as above $|\omega_j|\leq (2 \pi)^{\frac{n}{2}}
(\prod_k j_k^{\frac{1}{4}}) e^{-\frac{\delta r^2}{2}}\leq c
r^{\frac{n}{2}} e^{-\frac{\delta r^2}{2}}.$ This holds true for
all $\omega_j, \ j$ in terms of $r$. Specializing our work for
large $r$, we have that $\forall \varepsilon>0, \ \exists R,$ such
that $\forall r>R, |\omega_j| \leq
e^{-\frac{1}{2}(\delta-\varepsilon) r^2}$. Note that the factor of
$\varepsilon$ is used to compensate for the $\sqrt{2 \pi}
j_k^{\frac{1}{4}}$ terms. The probability of which may be
estimated using Lemma \ref{Gauss} as:
\par $$\nu(\{\omega: |\omega_{j}|\leq e^{-\frac{1}{2}(\delta-\varepsilon) r^2}\}) \leq e^{-(\delta-\varepsilon) r^2}.$$
\\ Hence $E_{\delta,r} $ is a small family of events as:\\
$\nu(\{ \omega : \log M_{r,\omega}' \leq (\frac{1}{2}-\delta)r^2
\})$ \par$\leq \nu( \{ \omega : |\omega_{j}| \leq
e^{-\haf(\delta-\varepsilon) r^2}, \ and \ j_k \ \in [(1- \sqrt
\frac{\delta}{n}) r^2, \  (1+ \sqrt \frac{\delta}{n})r^2 ]
\})$\par$ \leq (e^{-(\delta-\varepsilon)
r^2})^{(2\sqrt\frac{\delta}{n} r^{2})^n}= e^{-2
^n(1+o(\delta))\delta^{\frac{n+2}{2}} r^{2n+2} }= e^{-c_{1,\delta}
r^{2n+2}}$, using the independence of $\omega_{j}$. \\
$\displaystyle{ M_{r,\omega}\geq
M_{\frac{1}{\sqrt{n}}r,\omega}'\geq \haf r^2 -\delta r^2}$, except
for small events thus proving the lemma.

\end{proof}

 Results of this type can deceive one
into thinking of random holomorphic functions as $ e^{\haf z^2}$.
This absolutely is not the case, as they are weakly invariant with
respect to the isometries of the reduced Heisenberg group. In
particular, an analog of the previous theorem holds at any point
(whereas this will be false for $ e^{\haf z^2}$).

\begin{corr}\label{ValueEstimate} 
For all $\delta > 0$ and  
$z_0 \in \overline{B(0,r)}\backslash B(0,\frac{1}{2}r)$, there
exists $\zeta\in B(z_0,\delta r)$ s.t.
$$ \log |\psi_\omega(\zeta)| > \left(\frac{1}{2} - 3 \delta \right) |z_0 |^2$$
except for on a small family of events.
\end{corr}
\begin{proof}
By Lemma \ref{Anchises}:
$$\nu(\{\omega : \max_{z \in \partial B(0,r)}\log|\psi_\omega(z)|-\frac{1}{2}|z|^2\leq - \delta r^2 \})\leq e^{-c r^{2n+2} }$$
By Lemma \ref{Invariance}, we have that for $z_0 \in
B(0,r)\backslash B(0,\haf r), \ z\in B(z_0,\delta r)$:
$$\nu(\{\omega : \max_{z \in \partial B(0,\delta r)}\log|\psi_\omega(z-z_0)|-\frac{1}{2}|z-z_0|^2\leq - \delta (\delta r)^2 \})\leq e^{-cr^{2n+2} }$$
\\
Hence, $\exists \ z \in B(z_0,\delta r)$ s.t.
$\log|\psi_\omega(z-z_0)|-\frac{1}{2} |z-z_0|^2\geq - \delta
(\delta r)^2$, except for a small family of events.
\par By hypothesis, $|z_0|\in [\frac{1}{2}r, r)$, hence
$|z-z_0| \leq \delta r \leq \frac{1}{4}r=
\frac{r}{2}\frac{1}{2}\leq \frac{1}{2}|z_0|$\par

Hence, $|z_0 - z|^2 \geq |z_0|^2 -\delta r^2 \geq |z_0|^2
(1-2\delta)$\par Without loss of generality assume that $\delta<
\frac{1}{4}$.\par
\begin{tabular} {cl} $\log |\psi_\omega(z-z_0)|$ & $\geq
\frac{1}{2} |z-z_0|^2-\delta^3 r^2\geq |z_0|^2 \frac{1}{2}
(1-2\delta)^2 -4 \delta^3 |z_0|^2$ \\ & $\geq \frac{1}{2} |z_0|^2
-2\delta |z_0|^2- \frac{1}{4}\delta |z_0|^2$ \\& $ \geq
\frac{1}{2} |z_0|^2- 3 \delta |z_0|^2 $
\end{tabular}\\
And, setting $\zeta = z-z_0$ this is what we set out to prove.
\end{proof}

Using that $\log \displaystyle{ \max_{B(0,r)}} |\psi_\omega|$ is
an increasing function in terms of r, we have the following
corollary:

\begin{corr}\label{GrowthOfLogMax} For all $\delta > 0$

\begin{tabular}{rl}
a)&$\displaystyle{Prob\left( \left\{ \omega : \ \lim_{r\rightarrow
\infty} \frac{(\log \max_{z\in B(0,r)} |\psi_\omega (z)|)-\haf
r^2}{r^2} \notin [-\delta, \delta] \right\}\right) =0}$\\
b)&$\displaystyle{Prob\left( \left\{ \omega : \ \lim_{r\rightarrow
\infty} \frac{(\log \max_{z\in B(0,r)} |\psi_\omega (z)|)-\haf
r^2}{r^2} \neq 0 \right\}\right) =0}$\\
\end{tabular}\\
\end{corr}

This corollary as well as corollary \ref{FiniteOrder2} have
already been proven by more direct methods,
\cite{ShiffmanZelditchVarOfZeros}.

\begin{proof} Part b follows immediately from part a, which we now prove:\par
Let $E_{\delta, R} =\{\omega: \frac{\log \max_{B(0,R)}
\psi_\omega(z) - \haf R^2}{R^2} \notin [-\delta, \delta ]\}$\par
Let $R_m= r + \delta (m+1) r, \ r>0$\par
Let $s_m \in [R_{m-1}, R_m ].$\\
Claim: $\forall m> M_\delta, \ \forall s_m, \ E_{\delta, s_m}
\subset E_{\frac{1}{3}\delta, R_m} \bigcup E_{\frac{1}{3}\delta,
R_{m-1}}$\par Let $M_\delta = \max \{ M_{1,\delta}, \
M_{2,\delta}\}$, which may be specifically determined. \\
Case i: for $\omega\in E_{\delta, s_m}, \ \log
\displaystyle{\max_{B(0,s_m)}} \psi_\omega \geq \haf s_m^2
+\delta s_m^2$\\

\begin{tabular}{rl} $\log \displaystyle{ \max_{B(0,R_m)}} |\psi_\omega| $&$\geq \haf s_m^2
+\delta s_m^2$,\\
&$\geq \haf (1+ m\delta )^2 r ^2 +\delta (1+ m\delta )^2 r^2$
\\ & $>\haf R_m^2 +\frac{1}{3} \delta R_m^2, \ \forall m> M_{1,\delta} $
\end{tabular}
\par
Therefore, $\omega\in E_{\frac{\delta}{3}, R_m} $\\ Case ii: for
$\omega\in E_{\delta, s_m}, \ \log \displaystyle{\max_{B(0,s_m)}}
\psi_\omega \leq \haf s_m^2
-\delta s_m^2$\\

\begin{tabular}{rl} $\log \displaystyle{ \max_{B(0,R_{m-1})}} |\psi_\omega| $&$\leq \haf s_m^2
- \delta s_m^2$\\
&$\leq \haf (1+ (m-1) \delta )^2 r ^2 - \delta (1+ m\delta )^2
r^2$
\\ & $\leq \haf R_{m-1}^2 -\frac{1}{3} \delta R_{m-1}^2, \ \forall m> M_{2,\delta} $
\end{tabular}
\par
Therefore, $\omega\in E_{\frac{\delta}{3}, R_{m-1}},$\par Hence,
$\forall m> M_\delta$ and $ \forall s \in [R_{m-1}, R_m], \
E_{\delta, s}\subset E_{\frac{1}{3}\delta, R_{m-1}} \cup
E_{\frac{1}{3}\delta, R_m} $

Hence, $Prob (\bigcup_{s\in [R_{m-1}, R_m]} E_{\delta, s}) \leq 2
e^{-c_\delta r^{2n+2} m^{2n+2}}$, and\par $\Sum_{m\in \N} Prob
(\bigcup_{s\in[R_{m-1}, R_{m}]} E_{\delta, s})
 = \Sum_{m\in \N} e^{-c_\delta
m^{2n+2}}<\infty$, and the result follows.
\end{proof}

\newsection{The Second main lemma}\par

Essentially to prove the main theorem that we are working towards
we need only one more interesting lemma, Lemma \ref{Athena}, in
which we will give an estimate for $\int \log |\psi_\omega|$. This
will be proved first by obtaining a crude estimate for $\int |\log
|\psi_\omega||$, except for a small family of events, and then by
proving facts about the Poisson Kernel, which will allow me to
approximate using Riemann integration the first integral with
values of $\log|\psi_\omega(z)|$ at a number of fairly evenly
spaced points.
\par

In order to establish notation I state the following standard
result:
\begin{proposition}
\label{Kundun} For $\zeta \in
B(0,r)$, h a harmonic function\\
$$h(\zeta) = \int_{\partial B(0,r)} P_r(\zeta,z) h(z)
d\sigma_r (z) $$ where $d\sigma_r$ is the Haar measure of the
sphere $ S_r= \partial B(0,r)$ and $P_r$ is the Poisson kernel for
$B(0,r)$.
\par
\end{proposition}
A proof of this can be found in many standard text books,
\cite{Krantz}. It is convenient to normalize $\sigma_r$ so that
$\sigma_r(S_r)=1$. For this normalization, the Poisson Kernel is:
$$P_r(\zeta,z)= r^{2n-2}\frac{(r^2-|\zeta|^2)}{|\zeta-z|^{2n}}
$$

\begin{lemma} \label{oakley} For all $r>R_n, \ \int_{\partial B(0,r)} |\log(|\psi_\omega|)| d \sigma_r(z) \leq
(3^{2n}+1) r^2$ except for a small family of events.
\end{lemma}
\begin{proof}
By Lemma \ref{Anchises}, with the exception of a small family of
events, there exists $\zeta_0 \in
\partial B(0,\frac{1}{2}r)$ such that $\log(|\psi_\omega (\zeta_0)|) >
0$.\par Combining this with Proposition \ref{Kundun},
$$\int_{\partial B(0,r)} P_r(\zeta_0, z) \log(|\psi_\omega(z)|) d
\sigma_r (z)\geq \log(|\psi(\zeta_0)|)\geq 0.$$ Hence,
$$\int_{\partial (B(0,r))} P_r(\zeta_0, z) \log^-(|\psi_\omega(z)|)
\leq \int_{\partial (B(0,r))} P_r(\zeta_0, z)
\log^+(|\psi_\omega(z)|)$$
\par
Since $\zeta \in \partial B(0,\haf r)$ and $z \in
\partial
B(0,r)$, we have:  $\haf r \leq |z- \zeta | \leq \frac{3}{2} r$.\\
Hence by using the formula for the Poisson Kernel,
$$\frac{1}{3}\left(\frac{2}{3}\right)^{2n-2} \leq P_r(\zeta,
z)\leq (2)^{2n-2} 3$$ Therefore, $ \int_{\partial
B(0,r)}\log^+(|\psi_\omega(z)|) d \sigma_r (z)\leq \log M_r \leq
(\frac{1}{2}+ \delta) r^2\leq r^2$, except for a small family of
events, by
Lemma \ref{Anchises}.\\
\begin{tabular}{rl}
$\int_{\partial (B(0,r))} P(\zeta_0, z)\log^+(|\psi_\omega(z)|)$ &
$\leq \sigma_r(S_r) \log (M_r) 3 (2)^{2n-2} $\\ &$\leq 3 (2)^{2n-2} r^2  $\\
$\int_{\partial (B(0,r))} \log^-(|\psi_\omega(z)|)d \sigma_r(z)$ & $\leq \frac{1}{\min_{z} P(\zeta_0, z)} \int_{\partial (B(0,r))} P(\zeta_0, z)\log^+(|\psi_\omega(z)|) $\\
&
$\leq 3\left(\frac{3}{2}\right)^{2n-2} \int_{\partial (B(0,r))} P(\zeta_0, z)\log^+(|\psi_\omega(z)|) $\\
& $\leq
 9 \left(\frac{3}{2}\right)^{2n-2}(2)^{2n-2} r^2 $ \\ & $\leq
3^{2n} r^2 $
\end{tabular}\\
And, the result follows immediately.

 \end{proof}

As we are already able to approximate $\log |\psi_\omega(z)|$ at
any finite number of points in order to use Reimann integration to
prove Lemma \ref{Athena} we will need to be able to choose
"evenly" spaced points on the sphere, as chosen according to the
next proposition:\par

\begin{lemma}\label{Sphere} (A partition of a Sphere)\\
If $(2n) m^{2n-1} = N $ \\ then $S_r^{2n}\subseteq \R^{2n}$ can be
"divided" into measurable sets $\{I^r_1, I^r_2, \ldots, I^r_N\}$
such that:\par 1) $\bigcup_j I^r_j = S_r$
\par 2) $\forall j \neq k, \ I^r_j \bigcap I^r_k = \emptyset$, \par 3) $diam(I^r_j)\leq
\frac{\sqrt{2n-1}}{m}r=
  \frac{c_n}{N^{\frac{1}{2n-1}}}r$
\end{lemma}
\begin{proof}
Surround $S_r$ with 2n pieces of planes: $P_{+,1}, P_{+,2},\ldots,
P_{+,n}, P_{-,1},\ldots P_{-,n}$, where $$P_{+,j}= \{ x\in
\R^{2n+1} : \ ||x||_{L^\infty}=r, \ x_j=r \}
$$
$$P_{-,j}= \{ x\in
\R^{2n+1} : \ ||x||_{L^\infty}=r, \ x_j=-r \}
$$
Subdivide each piece into $m^{2n-1}$ even $2n-1$ cubes, in the
usual way, and denote these sets $R_1, \ldots, \ R_N$.\par Let
$I^r_j= \{ x\in S_r: \ \lambda x \in R_n, \ \lambda > 0 \}$\par By
design, $\lambda\geq 1$ and $x, \ y \in I_j \Rightarrow \ d(x,y)<
d(\lambda_1 x, \lambda_2 y)\leq \frac{2}{m} r= diam(R_j) $. These
sets can be redesigned to get that $I^r_j \bigcap I^r_k =
\emptyset$ , $j \neq k$ by carefully defining $R_i$ so that $R_j
\bigcap R_k= \emptyset$.
\end{proof}

The following elementary result is less well known then others and
will be very useful in proving Lemma \ref{Athena}. Note this
integration is with respect to w, which is not the same variable
of integration that is used in Proposition \ref{Kundun}. This is
done because the goal of this section, Lemma \ref{Athena}, is to
estimate a surface integral, which corresponds to integration with
respect to the first variable.

\begin{lemma} \label{poison}For $\kappa<1$
$$\int_{w \in S_{\kappa r}^n} P_r(w,z) d\sigma_{\kappa r}(w)= 1 $$
\end{lemma}

\begin{proof}
\par $P_r(w,z)= r^{2n-2}
\frac{r^2-|w|^2}{|z-w|^{2n}}$, $z \in
\partial B(0,r)$ \par If $w\in  S_{\kappa r}^{2n} \subseteq
\R^{2n}$, then the poisson Kernel can be rewritten as a function
of $|z-w|$, and as such $\forall \Upsilon \in U_n (\R^n), \
P_r(\Upsilon w, \Upsilon z )=P_r(w,z)$\par Let $f(z)= \int_{w \in
S_{\kappa r}^n} P_r(w,z) d\sigma_{\kappa r}(w)$\par

\begin{tabular}{cl} $f(z)$&
$=\int_{w \in S_{\kappa r}^n} P_r(w,z) d\sigma_{\kappa r}(w)$ \\
& $= \int_{w \in  S_{\kappa r}^n} P_r(\Upsilon w,\Upsilon z)
d\sigma_{\kappa r}(w)$, by the above work.\\ & $= \int_{w \in
S_{\kappa r}^n} P_r(\Upsilon w,\Upsilon z) d\sigma_{\kappa
r}(\Upsilon w)$, as $d\sigma_{\kappa r}$ is invariant under
rotations. \\ & $= \int_{w \in
S_{\kappa r}^n} P_r (w, \Upsilon z) d \sigma_{\kappa r}(w)$, by a change of coordinates.\\
& $= f(\Upsilon z)
$\end{tabular} \\ Hence $f(z)= c, \ \forall z \in S_r^n $\\
By switching the order of integration we compute that: $$1=
\int_{w \in S_{\kappa r}^n} \int_{z \in S_r^n}
P_r(w,z)d\sigma_r(z) d\sigma_{\kappa r}(w)= c$$
\end{proof}
Now we are able to prove our final lemma.

\begin{lemma}\label{Athena} For all $\Delta > 0,$\par
$ \left\{\omega : \frac{1}{r^2} \displaystyle\int_{z\in\partial
B(0,r) }  \log |\psi_\omega | d \sigma_r(z)  \leq \frac{1}{2} -
\Delta \right\}$ is a small family of events.\\
\end{lemma}
\begin{proof}
It suffices to prove the result for small $\Delta$. Let
$\Delta>0$. Let \newline$a_n = \frac{1}{2(2n +2) (2n-1)}$. Set
$\displaystyle{\delta=\left(\frac{1}{ \lambda}
\Delta\right)^{(\frac{1}{a_n})}<\frac{1}{6}}, \ \lambda>0$ to be
determined later. Choose $m \in \N$ such that writing $N = (2n)
m^{2n-1}, \ \frac{1}{N}\leq \delta$.  Let
$\kappa=1-\delta^{a_n}$.\par Choose $I^{\kappa r}_j$ measurable
subsets of $ S_{\kappa r}$ as in Proposition \ref{Sphere}. \\In
particular:\par
\begin{tabular}{cl} 1) & $S_{\kappa r} = \cup I^{\kappa r}_j$, a disjoint union.\\
2) & $\Sum \sigma_r(I^{\kappa r}_j)=1 $\\ 3) & $diam(I^{\kappa
r}_j)\leq   \frac{c_n}{N^{\frac{1}{2n-1}}} \kappa r\leq c
\delta^{\frac{1}{2n-1}}r$
\end{tabular}\\
Let $\sigma_j= \sigma_{\kappa r}(I^{\kappa r}_j)$, which does not
depend on r.\par For all $j$ fix a point $x_j \in I_j^{\kappa r}$

\par By Lemma \ref{ValueEstimate}, $\exists \zeta_j \in B(x_j,\delta r)$ such that
$$\log (|\psi_\omega(\zeta_j)|) > \left(\frac{1}{2} - 3 \delta \right)|x_j|^2=
\left(\frac{1}{2} - 3 \delta\right)\kappa^2 r^2 $$ Except, of
course, on N different small families of events (the union of
which remains a small family of events).

\begin{tabular}{cl}

\hspace{.1in} $\displaystyle{(\frac{1}{2}-3 \delta)}$ &
$\displaystyle{(1-\delta^{a_n})^2 r^2 \leq
\Sum_{j=1}^N \sigma_j \log(|\psi_\omega (\zeta_j)|)}$\\
& $\displaystyle{\leq \int_{\partial B(0,r)}  \left( \Sum_j
\sigma_j P_r(\zeta_j, z) \log(|\psi_\omega (z)|) d \sigma_r
(z)\right) }$
\\
& $\displaystyle{= \int_{\partial (B(0,r))} \left(\Sum_j \sigma_j
(P_r(\zeta_j, z) -1) \right) \log(|\psi_\omega(z)|) d \sigma_r (z)}$\\
&$\displaystyle{ \ \ \ \  + \int_{\partial (B(0,r))}
\log(|\psi_\omega (z)|) d \sigma_r (z)}$
\\
\end{tabular}\\
Hence,
\\\begin{tabular}{rl} $\int_{\partial B(0,r)} $& $ \log(|\psi_\omega|) d \sigma_r$\\&  $\geq
(\frac{1}{2}-3 \delta) (1-\delta^{a_n})^2 r^2 -\int
|\log|\psi_\omega|| d \sigma_r \cdot \max_z |\sum_j \sigma_j
(P_r(\zeta_j, z) -1) | $
\\ & $\geq (\frac{1}{2}-3 \delta) (1-\delta^{a_n}) r^2 -
(3^{2n}+1)r^2\cdot C_n \delta^{\frac{1}{2(2n-1)}} \geq \haf r^2 - \lambda \delta^{a_n} r^2$\end{tabular}\\
by Lemmas \ref{oakley} and the following claim. After proving this
claim, the result will follow.\\
Claim: $\displaystyle{ \max_{z\in
\partial (B(0,r))} \left|  \sum_j \sigma_j (P_r(\zeta_j,z)
-1 )\right| \leq C_n \delta^{\frac{1}{2 (2n-1)}}}$\\ Proof of
claim: $\forall z\in \partial B(0,r), \ \int_{\zeta \in
\partial B(0,\kappa r)} P_r(\zeta,z) d \sigma_{\kappa r} (\zeta) = 1$, by Lemma
\ref{poison}.
\par
Hence, $ 1 = \sum_{j=1}^{j=N} \sigma_j P_r(\zeta_j,z) +
\sum_{j=1}^{j=N} \int_{\zeta\in I^{\kappa r}_j}
(P_r(\zeta,z)-P_r(\zeta_j,z)) d\sigma_{\kappa r}(\zeta)$\par
\begin{tabular}{cl}
And, $ |\sum_{j=1}^{j=N}$&$ \sigma_j
(P_r(\zeta_j,z)-1)|=|\sum_{j=1}^{j=N} \int_{\zeta\in I^{\kappa
r}_j}
(P_r(\zeta,z)-P_r(\zeta_j,z)) d\sigma_{\kappa r}(\zeta)|$\\
&$ \displaystyle{\leq \max_{j, \ \zeta \in I^{\kappa r}_j}|\zeta -
\zeta_j| \cdot \max_{w \in B(0, (\kappa + \delta) r) \backslash
B(0, (\kappa - \delta) r)} \left|\frac{\partial P_r(w,z)}{\partial
w} \right| }$

\end{tabular}\par





$\displaystyle{\frac{\partial P_r( w, z )}{\partial w} = -r^{2n-2}
\frac{\overline{w} |z-w|^2 + (r^2-|w|^2) n
 (\overline{z}-\overline{w})
}{|z-w|^{2n+2}}}$\par As $|z|=r$, and $|w|=(1-\varepsilon) r\in
[(\kappa - \delta)r, (\kappa + \delta)r]$\par
\par $\left|\frac{\partial P_r( w, z )}{\partial w}\right|
\leq \frac{2+ 4 \varepsilon n}{r \varepsilon^{2n+2}}\leq
\frac{c_n}{r \varepsilon^{2n+2}}=\frac{c_n}{r }\delta^{
-\frac{1}{2(2n-1)}}$

And, $\max_{\zeta} |\zeta-\zeta_j| \leq diam(I_j)+ \delta r \leq c
\delta^{\frac{1}{2n-1}} r+\delta r\leq c' r \delta
^{\frac{1}{2n-1}}$
\par
Therefore: $| \sum_{j=1}^{j=N} \sigma_j ( P_r(\zeta_j, z)-1)| \leq
C \delta^{\frac{1}{2n-1}} \cdot \delta^{-\frac{1}{2(2n-1)}} =C
\delta^{\frac{1}{2 (2n-1)}}$\\ Proving the claim and the
lemma.\par

\end{proof}


This lemma gives an alternate proof for the growth rate of the
characteristic function. Let $T(f,r)= \int_{S_r} \log^+ |f(z)|
d\sigma_r(z) $, the Nevanlina characteristic function. As $(
\int_{S_r} \log |\psi_\omega| d\sigma_r)$ is increasing the proof
of Corollary \ref{GrowthOfLogMax} can be used in conjunction with
Lemma \ref{Athena} to prove that $\psi_\omega(z)$ is a.s. finite
order 2.

\begin{corr} \label{FiniteOrder2} For all $\delta \in (0,\frac{1}{3}]$\par
\begin{tabular}{rl}
a)&$\displaystyle{Prob\left( \left\{ \omega : \ \lim_{r\rightarrow
\infty} \frac{( \int_{S_r} \log |\psi_\omega| d\sigma_r)-\haf
r^2}{r^2} \notin [-\delta, \delta] \right\}\right) =0}$\\
b)&$\displaystyle{Prob\left( \left\{ \omega : \ \lim_{r\rightarrow
\infty} \frac{( \int_{S_r} \log |\psi_\omega| d\sigma_r)-\haf
r^2}{r^2} \neq 0 \right\}\right) =0}$\\
c)&$\displaystyle{Prob\left( \left\{ \omega : \ \lim_{r\rightarrow
\infty} \frac{T(\psi_\omega, r)-\haf
r^2}{r^2} \neq 0 \right\}\right) =0}$\\
\end{tabular}\\
\end{corr}

\newsection{Proof of Main results}\par

We will now be able to put the pieces together to estimate the
number of zeroes in a large ball for a random holomorphic function
$\psi_\omega(z)$. Further, This will help us to compute the hole
probability.

\begin{deff} For $f\in \holo (B(0,r)), \ B(0,r) \subset \C^n$, the
unintegrated counting function,\\
$n_{f}(r):= \int_{B(0,t)\bigcap Z_{f}} (\frac{i}{2 \pi} \partial
\overline{\partial} \log |z|^2)^{n-1}= \int_{B(0,t)} (\frac{i}{2
\pi} \partial \overline{\partial} \log |z|^2)^{n-1} \wedge
\frac{i}{2 \pi} \partial \overline{\partial} \log |f| $
\end{deff}
\par The equivalence of these two definitions follows by the
Poincare-Lelong formula. The above form ($(\frac{i}{2 \pi}
\partial \overline{\partial} \log |z|^2)^{n-1} $) gives a
projective volume, with which it is more convenient to measure the
zero set of a random function. The Euclidean volume may be
recovered as $\int_{B(0,t)\bigcap Z_{f}} (\frac{i}{2 \pi}
\partial \overline{\partial} \log |z|^2)^{n-1}=
\int_{B(0,t)\bigcap Z_{f}} (\frac{i}{2 \pi t^{2}} \partial
\overline{\partial} |z|^2)^{n-1}$.\par

\begin{lemma}\label{Shabat}If $u \in L^1(\overline B _r), \
and \ \partial \overline{\partial}u$ is a measure, then
$$\int_{t=r\neq 0}^{t=R} \frac{dt}{t} \int_{B_t} \frac{i}{2 \pi} \partial
\overline{\partial} u \wedge (\frac{i}{2 \pi} \partial
\overline{\partial} \log |z|^2)^{m-1}= \frac{1}{2} \int_{S_R}u
d\sigma_R - \frac{1}{2} \int_{S_r}u d\sigma_r
$$
\end{lemma}
\par A proof of this result is available in the literature, \cite{ShiffmanEquidistTheort}.\par When applying
this to random functions, my previous estimates of the surface
integral will turn out to be extremely valuable. \\
{\bf Theorem \ref{Main}} For all $\delta>0,$
$$F_r:=\left\{\omega: \left|n_{\psi_\omega}(r) -\frac{1}{2}r^2 \right| \geq \delta r^2 \right\} \ is \ a \ small \ family \ of \ events.$$

\begin{proof} It suffices to prove the result for small $\delta$.\par We will start by estimating that:
$$\nu \left( \left\{ \omega: \frac{n_{\psi_\omega}(r)}{r^2}  \geq \frac{1}{2} + \delta \right\}\right)
\leq e^{-c_\delta r^{2n+2}}$$

$n_{\psi_\omega}(r) \log(\kappa) \leq \int_{t=r}^{t=\kappa r}
n_{\psi_\omega}(t)\frac{dt}{t} \leq n_{\psi_\omega}(\kappa r)
\log(\kappa)$, as n(r) is increasing.\par let
$\kappa = 1+ \sqrt{\delta} $. Except for a small family of events, we have: \\

\begin{tabular}{cl} $\displaystyle{n_{\psi_\omega}(r) \log(\kappa)}$ & $\displaystyle{\leq
\int_{t=r}^{t=\kappa r} n_{\psi_\omega}(t)\frac{dt}{t} }$ \\ &
$\displaystyle{ = \int_{t=r}^{t= \kappa r} \int_{B(0,t)}
\frac{i}{2 \pi} \partial \overline{\partial} \log|\psi_\omega(z)|
\wedge \left(\frac{i}{2 \pi} \partial \overline{\partial}
\log|z|^2\right)^{n-1} \frac{dt}{t} }$\\
&$\displaystyle{ = \frac{1}{2} \int_{S_{\kappa r}} \log
|\psi_\omega (z)| d\sigma- \frac{1}{2} \int_{S_r} \log
|\psi_\omega (z)| d\sigma }$, by Lemma \ref{Shabat}.\\ &
$\displaystyle{ \leq
\frac{1}{2}\left(\left(\frac{1}{2}+\delta\right) \kappa^2 r^2 -
\int_{S_r} \log |\psi_\omega (z)| d\sigma\right) }$, by Lemma
\ref{Anchises}. \\ & $\displaystyle{ \leq
\frac{1}{2}\left(\left(\frac{1}{2}+\delta\right) r^2 \kappa^2 -
\left(\frac{1}{2} - \delta\right)r^2\right)} $, by Lemma
\ref{Athena}.
\end{tabular}

\begin{tabular}{rl}$\displaystyle{2 \frac{n_{\psi_\omega}(r)}{r^2}}$&$\displaystyle{\leq
\frac{1}{\log (\kappa)} \left(\kappa^2 \left(\frac{1}{2}
+\delta\right)-\left(\frac{1}{2}-\delta \right) \right) }$\\&
$\displaystyle{= \frac{\kappa^2 -1}{2 \log(\kappa)}+ \delta
\frac{\kappa^2+1}{\log(\kappa)} \leq 1 + c \sqrt{\delta}}$.\end{tabular}\\

This proves the probability estimate when the unintegrated
counting function is significantly larger then expected.\par In
order to prove the other probability estimate:
$$\nu\left( \left\{\omega: \frac{n_{\psi_\omega}(r)}{r^2} \leq \frac{1}{2} - \delta \right\} \right)
\leq e^{-c_\delta r^{2n+2}}$$

We start by using that: $\int_{t=\kappa^{-1} r}^{t= r}
n_{\psi_\omega}(t)\frac{dt}{t} \leq n_{\psi_\omega}( r)
\log(\kappa)$. We then use that, except for a small family of
events, we have that:

\begin{tabular}{cl} $\displaystyle{n_{\psi_\omega}(r) \log(\kappa)}$ & $\displaystyle{\geq
\int_{t=\kappa^{-1} r}^{t= r} n_{\psi_\omega}(t)\frac{dt}{t}}$ \\
& $\displaystyle{ = \int_{t=\kappa^{-1} r}^{t=r} \int_{B(0,t)}
\frac{i}{2 \pi}
\partial \overline{\partial} \log|\psi_\omega(z)| \wedge
\left(\frac{i}{2 \pi} \partial \overline{\partial} \log|z|^2
\right)^{n-1}\frac{dt}{t}}$
\\ &$\displaystyle{= \frac{1}{2} \int_{S_r} \log |\psi_\omega (z)|
d\sigma- \frac{1}{2}
\int_{S_{\kappa^{-1}r}} \log |\psi_\omega (z)| d\sigma } $, by Lemma \ref{Shabat}.\\
& $\displaystyle{\geq \frac{1}{2}[(\frac{1}{2}-\delta) r^2 -
\int_{S_{\kappa^{-1}r}} \log |\psi_\omega (z)| d\sigma]}$, by
Lemma \ref{Athena}.
\\ & $\displaystyle{\geq \frac{1}{2}[(\frac{1}{2}-\delta) r^2 - (\frac{1}{2} +
\delta)r^2\kappa^{-2}]}$, by Lemma \ref{Anchises}.
\end{tabular}

\begin{tabular}{rl}$\displaystyle{2 \frac{n_{\psi_\omega}(r)}{r^2}}$&$\displaystyle{\geq \frac{1}{\log
(\kappa)} \left( \left(\frac{1}{2}
-\delta\right)-(\frac{1}{2}+\delta) \kappa^{-2} \right)} $\\
&$\displaystyle{= \frac{1-\kappa^{-2}}{2 \log(\kappa)}- \delta
\frac{1+\kappa^{-2}}{\log(\kappa)} \geq 1- 2 \sqrt \delta }$
\end{tabular}
\\ 

\end{proof}

Using this estimate for the typical measure of the zero set of a
random function we get an upper bound for the hole probability,
and putting this together with some elementary estimates we get an
accurate estimate for the order of the decay of the hole
probability:\\
{\bf Theorem \ref{Hole probability}} {\it If}
$$\psi_\omega (z_1, z_2 \ldots, z_n)= \Sum_j \omega_j
\frac{z_1^{j_1} z_2^{j_2} \ldots z_n^{j_n}}{\sqrt{j_1! \cdot
j_n!}},$$ {\it where $\omega_j$ are independent identically
distributed complex Gaussian random variables, and } $$Hole_r=\{
\omega:\forall z \in B(0,r), \ \psi_\omega (z)\neq 0 \},$$ {\it
then there exists} $c_1, c_2
>0$ {\it such that for all} $r>R_n$
$$ e^{-c_2 r^{2n+2}} \leq Prob (Hole_r)\leq e^{-c_1 r^{2n+2}}$$


\begin{proof}  The upper estimate follows by the previous
theorem, as if there is a hole then $n_{\psi_\omega}(r)=0$, and
this can only occur on a small family of events.\par Therefore it
suffices to show that the hole probability is bigger than a small
set.\par
 Let $\Omega_r$ be the event
where:\\
$i) \ |\omega_{0}| \geq E_n + 1$,\\
$ii) \ |\omega_{j}|\leq e^{-(1+ \frac{n}{2})r^2}, \ \forall j:
1\leq
|j|\leq \lceil 24 n r^2 \rceil= \lceil(n \cdot 2 \cdot  12)r^2 \rceil $\\
$iii) \ |\omega_{j}|\leq 2^{\frac{|j|}{2}},  \ |j| > \lceil 24 n r^2 \rceil \geq 24 n r^2 $\\
$\nu( \{\omega | \ |\omega_{j}|\leq e^{-(1+ \frac{n}{2})r^2} \})
\geq \frac{1}{2}(e^{-(1+ \frac{n}{2})r^2})^2=\frac{1}{2} e^{-(2+
n)r^2} $, by Lemma \ref{Gauss} \par

$ \# \{j \in \ \N^n | 1\leq |j|\leq \lceil 24 n r^2 \rceil \} =
({\lceil 24 n r^2 \rceil + n \choose n} ) \approx  c r^{2n} $\par
Hence, $\nu(\Omega_r)\geq C (e^{-c_n r^{2n+2}})$, by independence and Lemma \ref{Gauss}. Therefore $\Omega_r$ contains a small family of events, and it now suffices to show that for $\omega \in \Omega_r, \ \psi_\omega$ has a hole in $B(0,r)$.  \\
$f(z) \geq |\omega_{0}| - \Sum_{|j|=1}^{|j| \leq \lceil 24 n r^2
\rceil} |\omega_{j}| \frac{r^{|j|}}{\sqrt{j!
 }} - \Sum_{|j| > \lceil 24 n r^2 \rceil} |\omega_{j}| \frac{r^{|j|}}{\sqrt{j!}}= |\omega_{0}|- \sum^1 - \sum^2$\\

\begin{tabular}{cl}
 $\Sum^1$ & $\leq e^{-(1+ \frac{n}{2})r^2} \Sum_{|j| = 1}^{|j| \leq \lceil 24 n r^2 \rceil} \frac{r^{|j|}}{\sqrt{j!}}
 $ \\
 \vspace{.1in}
 &$\leq e^{-(1+ \frac{n}{2})r^2} \sqrt{(24 n r^2 +1)^{n}} \sqrt{(e^{r^n})}$,
 by Cauchy-Schwarz inequality.
\\
 \vspace{.1in}
 &$\leq C_n r^{n} e^{-r^2}\leq c e^{-0.9 r^2}<\frac{1}{2} $ for $r> R_n$\\
\end{tabular}\\

\begin{tabular}{cl}
 $\Sum^2$ & $\leq \Sum_{|j| > 24 n r^2} 2^{\frac{|j|}{2}} \left(\frac{|j|}{24 n}\right)^{\frac{|j|}{2}} \frac{1}{\sqrt{j!
}}$, as $r<\sqrt{\frac{|j|}{24 n}} $\\
  & $\leq c \Sum_{|j| > 24 n r^2} 2^{\frac{|j|}{2}} \left(\frac{|j|}{24 n}\right)^{\frac{|j|}{2}} \prod_{k=1}^{k=n} \left(\frac{e}{j_k}\right)^{\frac{j_k}{2}}$, by Sterling's formula \\
  & $=c\Sum_{|j| > 24 n r^2} \frac{(|j|)^{\frac{|j|}{2}}}{\left(\prod_{k=1}^{k=n} j_k^\frac{j_k}{2}\right) n^{\frac{|j|}{2}}} \left(\frac{e}{12}\right)^{\frac{|j|}{2}}
  $\\
& $\leq c \Sum_{|j|>1} \left(\frac{1}{4}\right)^{{\frac{|j|}{2}}}$, by Lemma \ref{CalcII}.\\
& $ \leq c \Sum_{l>1} \left(\frac{1}{2}\right)^{l} l^n\leq E_n $
\end{tabular}\\
Hence, $|\psi_\omega(z)| \geq E_n + 1- \Sum^1-\Sum^2\geq \frac
{1}{2} $
\end{proof}

\end{document}